\documentclass[12pt]{iopart}
\usepackage{amsfonts}
\usepackage{cite}
\usepackage{graphicx}
\usepackage{color}

\newcommand{\divv}{\mathrm{div}}
\newcommand{\curl}{\vec{\mathrm{curl}}}

\begin{document}

\title{Practical error estimates for sparse recovery in linear inverse
problems}
\author{Ignace Loris and Caroline Verhoeven}
\address{Mathematics Department, Vrije Universiteit Brussel, Pleinlaan 2, 1050 Brussel, Belgium}
\ead{igloris@vub.ac.be,cverhoev@vub.ac.be}

\begin{abstract}
The effectiveness of using model sparsity as a priori
information when solving linear inverse problems is studied. We
investigate the reconstruction quality of such a method in the
non-idealized case and compute some typical recovery errors
(depending on the sparsity of the desired solution, the number
of data, the noise level on the data, and various properties of
the measurement matrix); they are compared to known
theoretical bounds and illustrated on a magnetic tomography
example.
\end{abstract}



\maketitle

\section{Introduction}

Compressed sensing \cite{Donoh2006,Bruckstein.Donoho.ea2009} is
an area of signal recovery which has recently attracted a great
deal of attention thanks to its large potential for
applications; it enables reconstruction of sparse signals
$x_0\in\mathbb{R}^n$ with far fewer linear measurements than
traditionally required. In this note we investigate what
happens to reconstruction quality when theoretical conditions
on the measurement matrix are relaxed and replaced by (less
favorable) conditions that are encountered in realistic inverse
problems.

When a signal $x_0$ is known to be sparse (i.e. has few nonzero
components), a practical method to obtain accurate
reconstruction results from measurements $Kx_0$ is the
$\ell_1$ minimization
technique\cite{CaRoTa2006,Donoh2006,Chen.Donoho.ea1998}:
\begin{equation}
\tilde{x}=\arg\min_{Kx=Kx_0}\|x\|_1,\label{cs}
\end{equation}
($K$ is a known $m\times n$ operator). $\|x\|_1$ denotes the
$\ell_1$-norm of a vector: $\|x\|_1=\sum_i|x_i|$. Under certain
conditions it can be shown that the $\ell_1$ reconstruction
$\tilde x$ equals the unknown sparse input signal $x_0$ exactly
\cite{Candes2008}. Data obtained from real measurements are
more often than not corrupted by significant quantities of
noise; this means that $Kx_0$ remains unknown and instead
$y=Kx_0+\mathrm{noise}$ is acquired. In this case, the
ill-posed problem can be regularized by adding a sparsity
promoting $\ell_1$-norm penalty to a quadratic misfit
functional (assuming Gaussian noise) \cite{Daubechies2004b}. In
other words, the input signal can then be recovered from the
noisy data $y$ by minimizing the convex functional
\begin{equation}
\label{l1pen}
\bar{x}(\lambda)=\arg\min_x\|Kx-y\|^2+2\lambda\|x \|_1
\end{equation}
for a suitable choice of the penalty parameter $\lambda$.
$\|\cdot\|$ stands for the $\ell_2$-norm: $\|a\|^2=\sum_i
|a_i|^2$.
Contrary to the noiseless case, the recovered signal
(\ref{l1pen}) will not be exactly equal to the input signal
$x_0$.

The potential of the $\ell_1$ minimization (\ref{cs}) and
$\ell_1$ penalization method (\ref{l1pen}) for the
reconstruction of sparse signals has already been assessed
extensively, both from a theoretical point of view and with
numerical simulations
\cite{DoE:2003,Donoho2006,DoT:2008,Donoho.Elad.ea2006,ElB:2002,Fuc:2004}.
Nonetheless, most research focuses almost exclusively on
matrices which exhibit mutually incoherent columns or matrices
which (likely) satisfy the `Restricted Isometry Property' (RIP)
\cite{CaT:2005}, such as e.g. random matrices, structured
random matrices (e.g. random rows of a Fourier matrix) or
matrices composed of the union of certain bases. In most
practical situations it is difficult to verify if a matrix
satisfies the RIP. Here we set out to study the quality of
reconstruction when the measurement matrix $K$ does \emph{not}
fit this ideal category. One may e.g. have good a priori
control on the sparsity of the desired model (by the choice of
a suitable basis), but one may not have sufficient control on
the physical measurements to make the matrix satisfy the RIP.
Indeed, a single matrix column with small or zero norm will
destroy the RIP, but attempting sparse recovery may still be
appropriate or worthwhile in practice.

In this paper, we investigate the influence of the noise level
in the data and of the singular value spectrum of the
measurement matrix $K$ on the ability of the $\ell_1$ method
(\ref{l1pen}) to faithfully reconstruct a sparse input signal
$x_0$. We also discuss the behavior of the recovery as a
function of the sparsity of the input signal and as a function
of the indeterminacy of the system (number of rows with respect
to number of columns in $K$); we provide practical predictions
on the relative error of the recovered sparse signal with
respect to the sparse ground truth model. In order to
accomplish this, a large number of numerical simulations is
performed and the results are displayed using the diagrams
introduced in \cite{DTD:2006,DoT:2008}. We also provide an
illustration of the practical use of the method in a magnetic
tomography setting.

\section{Assessment method}
\label{sec:met}

We consider an input signal $x_0\in\mathbb{R}^n$ with $k$
non-zero coefficients, and a data vector $y\in\mathbb{R}^m$,
with $m\leq n$, such that $y=Kx_0+\eta$ for a $m\times n$
matrix $K$ and a noise vector $\eta$. We want to assess
the effectiveness of the $\ell_1$ penalization method
(\ref{l1pen}) for recovering $x_0$, using the knowledge of $K$
and of $y$, as a function of the spectrum of the matrix
$K$, of the noise level $\epsilon=\|\eta\|/\|Kx_0\|$, of the
number of data $m$ and of the sparsity $k$ of $x_0$.

The success of the $\ell_1$ method (\ref{l1pen}) for the
recovery of a sparse signal $x_0$ will depend on the
indeterminacy of the linear system and on the sparsity of the
input signal. A concise graphical representation of the success
rate of $\ell_1$-based compressed sensing was introduced in
\cite{DTD:2006,DoT:2008}. Likewise, we will use the parameters
$\delta=m/n$ and $\rho=k/m$ (not $k/n$) and perform a number of
experiments for various values of $\delta$ and $\rho$. More
precisely, we will use a cartesian grid in the
$\delta-\rho$-plane and, for each grid point, set up input
data, a matrix $K$ and a noise vector $\eta$, determine the
minimizer $\bar{x}(\lambda)$ and compare it to $x_0$. This
experiment is repeated several times over, for different $K$,
$x_0$ and $\eta$, but for fixed spectra and values of
$\epsilon$. Afterwards the whole experiment is repeated for
different spectra and different values of $\epsilon$. Do note
that the $\delta-\rho$-plane does not give a uniform
description of all possible inverse problems. It is
specifically tailored towards the evaluation of sparse
recovery.

In this work we use 40 equidistant points for $\delta$ and
$\rho$ with values ranging from 0.025 to 1. The matrices $K$
have a fixed number of columns ($n=800$). The experiment is
repeated 100 times in each grid point, for each spectrum and
each value of $\epsilon$. The same experiments with a larger
number of columns allow us to conclude that the result
presented in this paper do not depend on this quantity.

In case of noiseless data, the success rate of the recovery
strategy can be measured by simply computing the proportion of
successful reconstructions ($\tilde x=x_0$). In the case of
noisy data, one will \emph{never} have perfect reconstruction
($\bar x(\lambda)=x_0$), and therefore `good recovery' needs to
be defined in a different way. We measure the success of
recovery by calculating the mean of the relative reconstruction
error $e=\|\bar{x}(\lambda)-x_0\|/\|x_0\|$ over several trials.
It is this number that will be plotted as a
function of $\delta=m/n$ and $\rho=k/m$.

\begin{figure}
\begin{center}
\resizebox{\textwidth}{!}{\includegraphics{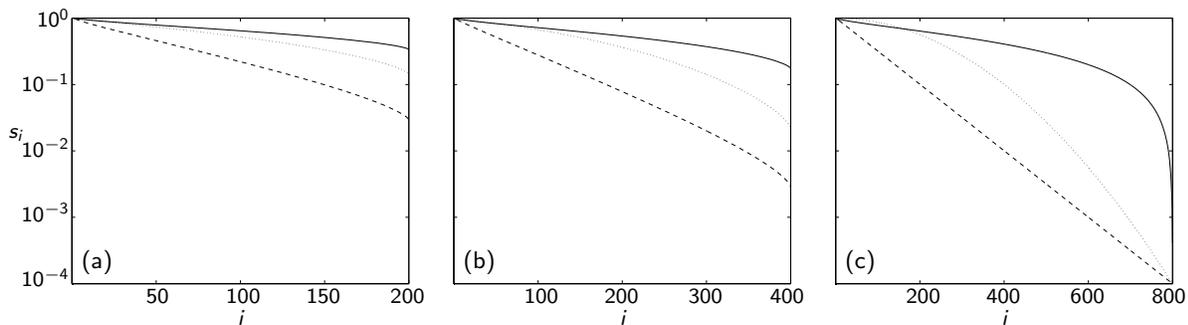}}
\end{center}
\caption{Mean over 100 examples of the normalized spectra of
matrices (see section \ref{sec:met}) of
the following type: random (solid), type 2 (dashed) and
type 3 (dotted), $\kappa=10^4$. The matrices are of size $200\times 800$
(a), $400\times 800$ (b) and $800\times 800$ (c).}
\label{fig:spect}
\end{figure}

When assessing the influence of the spectral properties of $K$
on the success of the method (\ref{l1pen}), we take into
account the change in spectrum of a matrix by the addition of
extra measurements (rows). We therefore choose the $m\times n$
matrices $K$ as submatrices of $n\times n$ matrices in the
following way: We draw an $n\times n$ matrix containing random
numbers from the Gaussian distribution with zero mean and unit
variance. Three different types of spectral behavior are then
considered. ``Type 1'' is obtained by taking the random matrix
itself. For the other types, we first calculate the singular
value decomposition and replace the singular values by
\begin{equation}
s_i=s_1\kappa^{(1-i)/(n-1)}\qquad i:1\dots n,
\label{spec:geo}
\end{equation}
(with $s_1\neq 0$) for ``type 2'', and by
\begin{equation}
s_i=s_1\kappa^{(1-i^2)/(n^2-1)}\qquad i:1\dots n,
\label{spec:gaus}
\end{equation}
for ``type 3''. Several values of $\kappa$
 will be chosen in the next section. The singular vectors
remain untouched. For all three types, the $m\times n$
matrices $K$ are then found by randomly selecting $m$
different rows from these $n\times n $ matrices. As seen in
figure \ref{fig:spect}, the spectrum changes with $m$. We
believe this corresponds to the behavior encountered in
real problems.

The position of the $k$ non-zero entries of the input vector
$x_0$, as well as their values,  were randomly
generated from a uniform distribution.

For a given matrix $K$ and a given sparse input signal $x_0$,
synthetic data are constructed by setting $y=Kx_0+\eta$ where
$\eta$ is a $m\times 1$ vector with entries taken from a
Gaussian distribution with zero mean. The noise level is
determined by the parameter $\epsilon=\|\eta\|/\|Kx_0\|$.
We will choose $\epsilon=0.02,0.05,0.10,0.20$ and $0.50$.

The penalty parameter $\lambda$ in (\ref{l1pen}) is chosen to
satisfy Morozov's discrepancy principle: $\|K\bar
x(\lambda)-y\|=\|\eta\|$. In other words, we will fit the data
up to the level of the noise. In practice, and for these
problems sizes, this can be achieved easily by using the
(non-iterative) Homotopy/LARS/lasso method
\cite{EfHJT2004,OsPrT2000,IMM2005-03897} for the solution of
(\ref{l1pen}).

\section{Results}
\label{sec:res}

\begin{figure}
\begin{center}
\resizebox{0.8\textwidth}{!}{\includegraphics{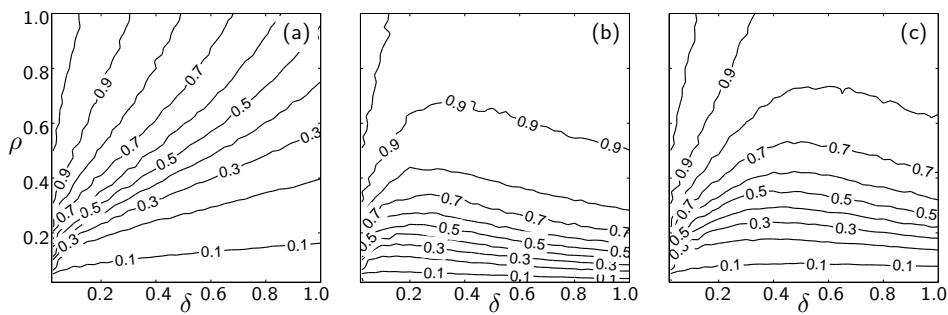}}
\caption{Mean, over $100$ trials, of the relative error
$e=\|\bar{x}(\lambda)-x_0\|/\|x_0\|$ as a function of $\delta$
and $\rho$ for a noise level of $10\%$. (a) Gaussian random
matrices, (b) matrices of type 2 and (c) matrices of type 3
with $\kappa=10^4$.}
\label{fig:recon}
\end{center}
\end{figure}

In order to assess the effect of the different variables on the
accuracy of the solution of (\ref{l1pen}), the mean of the
reconstruction error $e$ is plotted in the $\delta-\rho$-plane
for the three classes of matrices described in section
\ref{sec:met} with $\kappa=10^4$ and noise level of $10\%$ (see
figure~\ref{fig:recon}). The quality of the $\ell_1$
penalization reconstruction depends strongly on the spectrum of
the matrix $K$. In particular we see that good results for
random matrices must be considered as a best case scenario. The
results are less favorable for the two other types of matrices.

\begin{figure}
\begin{center}
\resizebox{\textwidth}{!}{\includegraphics{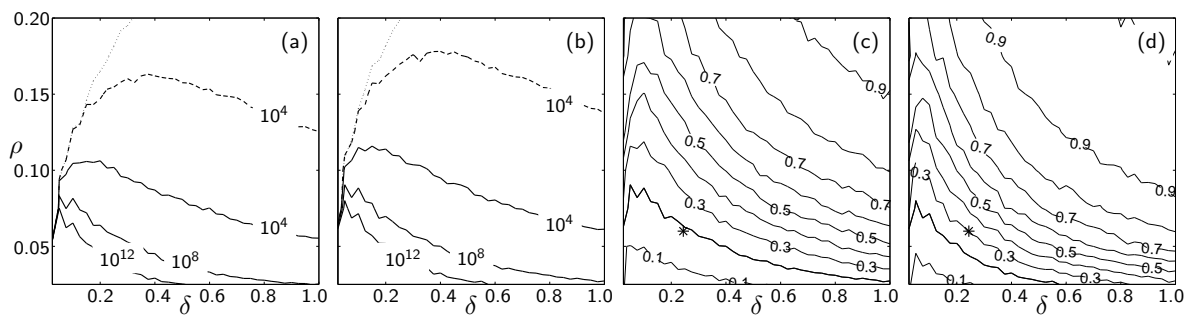}}
\end{center}
\caption{(a) and (b): Behavior in the $\delta-\rho$-plane of the mean
relative error $e_{2\epsilon}$ corresponding to two times the
noise level $\epsilon$ with respectively $\epsilon=2\%$ and $10\%$.
The reconstructions were done for the following matrices: Gaussian
random (dotted),type 2 (solid) and type 3 (dashed). The value of $\kappa$
is displayed on the curve. Panels (c) and (d): Mean relative errors
on reconstructions for matrices of type 2 with respectively $\kappa=10^8$
and $\kappa=10^{12}$ and $\epsilon=0.1$. The stars indicate
the position in the $\delta-\rho$-plane of the input model of an inverse problem in
magnetic tomography described in section \ref{sec:meg}.}
\label{fig:noisetest}
\end{figure}

Secondly, in order to combine results for different spectra of
$K$ in one plot, the behavior of the mean relative error
$e_{2\epsilon}=2\epsilon$ is studied. In figure
\ref{fig:noisetest} (a) and (b), the relative accuracy
$e_{2\epsilon}$ is plotted for various matrices $K$ and
$\epsilon=2\%$ and $10\%$ respectively. Only very sparse
vectors $x_0$ are shown in these plots ($0.025\leq\rho\leq
0.2$) as the results deteriorate rapidly for larger $\rho$.
These two $e_{2\epsilon}$ curves in the $\delta-\rho$-plane
(for the noise levels of $2\%$ and $10\%$) are quite similar in
character. We also simulated noise levels of $5\%$,
$20\%$ and $50\%$ and observed the same features.\\
The same sparse setting as figure 3(a,b) was used in figure 3
(c) and (d) where we show the reconstruction errors $e$ for
matrices of type 2 with $\kappa=10^8$ and $\kappa=10^{12}$
respectively, and $\epsilon=0.1$.

It is not surprising that the above results depend on the
spectrum of $K$. Any change in the spectrum of a matrix $K$
also induces a change in the condition number of the column
submatrices of $K$. As the restricted isometry property (RIP)
\cite{CaT:2005} implies that these submatrices are well
conditioned, our setting makes it more
difficult/impossible for a matrix $K$ to satisfy a RIP. As
argued below, even when the RIP is not satisfied, a sparse
recovery may still be feasible.

In \cite{Rau:2010} it is stated that when
\begin{equation}
(1-\delta_{2k})\|z\|\leq\|Kz\|\leq(1+\delta_{2k})\|z\|,\quad
\forall z: \|z\|_0\leq 2k,
\label{con:rip}
\end{equation}
and $\delta_{2k}\leq 2/(2+\sqrt{7/4})$ ($\|z\|_0$ is
the number of non-zero coefficients of $z$), any vector
$x_0$, with $\|x_0\|_0\leq k$ can be recovered by $\ell_1$
minimization with an error $\|\bar{x}-x_0\|\leq c\|\eta\|$. The
constant $c$ only depends on $\delta_{2k}$ and $\eta$ is the
noise vector as defined in section \ref{sec:met}. This
condition implies that the condition numbers ($\kappa_{2k}$) of
all the $2k$ column submatrices of $K$ are bounded from above
by:
\begin{equation}
\kappa_{2k}\leq \sqrt{\frac{10+\sqrt{7}}{2+\sqrt{7}}}\approx 1.6498,
\label{rip:con}
\end{equation}

As an illustration, we consider a matrix $K$ (of type 1, 2, 3
resp.) of size $200\times 800$ ($\delta=0.25$) and choose 10000
random 20-column submatrices of $K$. The average condition
number $\langle\kappa_{20}\rangle$ of those submatrices is
calculated: $\langle\kappa_{20}\rangle\approx 1.8$ for type 1
matrices, $\langle\kappa_{20}\rangle\approx 2.7$ and
$\langle\kappa_{20}\rangle\approx 2.0$ for matrices of type 2
and 3 respectively with $\kappa=10^4$. We see that
$\langle\kappa_{20}\rangle$ lies above the upper bound
(\ref{rip:con}). This means that for the three types of
matrices studied here, there exists some column submatrices
with 20 columns, which do not satisfy condition
(\ref{rip:con}). The matrix $K$ (even normalized) will
therefore not satisfy (\ref{con:rip}) for $k=10$ (which
corresponds to $\rho=0.05$) and therefore
\cite{Rau:2010} does not guarantee the accuracy of the
reconstruction of a vector with 10 (or more) nonzero
coefficients. However, as seen on figure \ref{fig:recon} and
\ref{fig:noisetest}, the $\ell_1$ penalization performs
reasonably well, even when the vector to be reconstructed is
less sparse than $\rho=0.05$. We can therefore conclude that
the bounds obtained from the RIP theory are quite restrictive.

The same conclusion can be drawn for the results obtained in
\cite{Donoho.Elad.ea2006}, where the behavior of the solution of
(\ref{l1pen}) with respect to the mutual
coherence of the matrix $K$ is studied. Inspired by
their analytical results, the authors only studied input
signals with a very small (maximum) number of nonzero
coefficients ($k=3$) for a matrix of size $128\times 256$. This
matrix is a concatenation of the identity matrix and the Hadamard
matrix, with the columns normalized to unit $\ell_2$ norm. For
that case we have calculated the value of the mean of
the condition number $\langle\kappa_{2k}\rangle$ over $10^6$
submatrices with 12 columns of $K$ (i.e. $k=6$ and $\rho=0.0469$)
and found $1.4216$, with a maximum of $1.7495$. From the present
simulations we therefore conclude that the $\ell_1$
penalization method can still yield accurate results for less
sparse solutions than imposed by the bound in
\cite{Donoho.Elad.ea2006}.

Most of the vectors $x_0$ studied here are less sparse than the
bounds found in the theoretical studies mentioned above.
Therefore, our simulations give an indication of the accuracy
of the results obtainable in practice with $\ell_1$
penalization for various types of matrices. They are of
practical use to predict the relative reconstruction error in
realistic inverse problem settings.

\section{Illustration}
\label{sec:meg}

We illustrate the previous results on the predicted accuracy of
sparse recovery by $\ell_1$ penalization with the help of an
inverse problem in magnetic tomography. Our toy problem aims to
reconstruct a 2D current distribution $\vec{J}$ on part of a
sphere from measurements of the (normal component of the)
magnetic field above this surface. The current distribution
$\vec{J}$ and the data $\vec{B}$ are assumed to be
linked by the formula \cite{Fornasier.Pitolli2008,VZP:1996}:
\begin{equation}
\vec{B}(\vec{r})=\frac{\mu_0}{4\pi}\int_V
\vec{J}(\vec{r}')\times
\frac{\vec{r}-\vec{r}'}{|\vec{r}-\vec{r}'|^3}
\,\mathrm{d}V'.\label{biotsavart}
\end{equation}
We also assume that the unknown currency distribution $\vec{J}$
is divergence-free: $\divv \vec{J}=0$.

In our toy problem the domain $V$ is part of a thin spherical
shell $0.089\leq r \leq 0.090$ (all units are SI units) which
we parametrize by coordinates $\xi$ and $\eta$:
\begin{equation}
x =\frac{r}{s} \tan\xi,\quad
y =\frac{r}{s} \tan\eta,\quad
z =\frac{r}{s}\quad\mathrm{and}\quad s=\sqrt{1+\tan^2\xi+\tan^2\eta}
\end{equation}
(with $-\pi/3\leq \xi,\eta \leq\pi/3$). The coordinate lines of
this parametrization are angularly equidistant great circles
\cite{Ronchi.Iacono.ea1996}. The shell covers just over one
quarter of the whole sphere. Furthermore we assume that the
current distribution has no radial component. The
divergence-free angular current distribution
$\vec{J}(\xi,\eta)$ is then parametrized by the field
$F(\xi,\eta)$: $\vec{J}\equiv \curl\times F(\xi,\eta)\vec{1}_r$
(the dependence on $r$ is not taken into account here).

The domain $V$ is divided into $64^2$ voxels and we choose an
input model $F^\mathrm{in}(\xi,\eta)$ that is sparse in the CDF
4-2 wavelet basis \cite{Cohen.Daubechies.ea1992}. The input
model has only $60$ nonzero coefficients in this basis (out of
a possible $4096$). Formula (\ref{biotsavart}) is used to
calculate the (normal component of the) magnetic field in 1000
randomly distributed points above the patch $-\pi/3\leq
\xi,\eta \leq\pi/3$ at $r=0.1$ (one centimeter above the
patch). 10\% Gaussian noise is added to these data. On this
$64\times64$ grid we therefore have 1000 data versus 4096
unknowns $F(\xi_i,\eta_i)$.

\begin{figure}\centering
\resizebox{\textwidth}{!}{\includegraphics{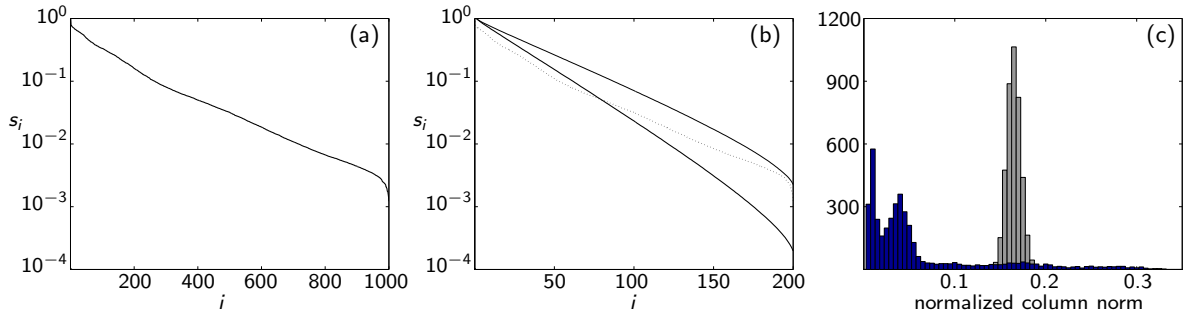}}
\caption{Left: The spectrum of the matrix used in the toy magnetic
tomography experiment. Center: The singular value
spectra of matrices of type 2 (solid) with $\kappa=10^8$ (top)
and $\kappa=10^{12}$ (bottom) and a rescaled spectrum of the matrix
of the toy problem (dotted).
Right: A histogram of columns norms of the toy problem matrix (blue)
and the matrix of (b) with $\kappa=10^{12}$ (grey) . They are both
normalized by their respective matrix spectral norms.}\label{megspectrum}
\end{figure}

The singular value spectrum of this $1000\times4096$ matrix (in
the wavelet basis) is plotted in figure~\ref{megspectrum} where
it is compared to the spectrum of the matrices used in sections
\ref{sec:met} and \ref{sec:res}. Its spectrum is quite similar
to the spectrum of a matrix of type 2 with $\kappa=10^{12}$.
Although the spectra are alike, it is important to point out
that some of the (wavelet) basis functions lie close to the
null space of the matrix. In other words, the norms of the
columns of the toy problem matrix are quite different from
those of the matrices used in the setting of the preceding
sections (see figure~\ref{megspectrum}). The average
reconstruction errors found in section~\ref{sec:res} should
therefore be understood as a lower bound: If one or more of the
basis functions (used to express the sparsity of the model)
lies more or less in the null space of the operator, the
reconstruction quality may be further reduced.

\begin{figure}\centering
\resizebox{\textwidth}{!}{\includegraphics{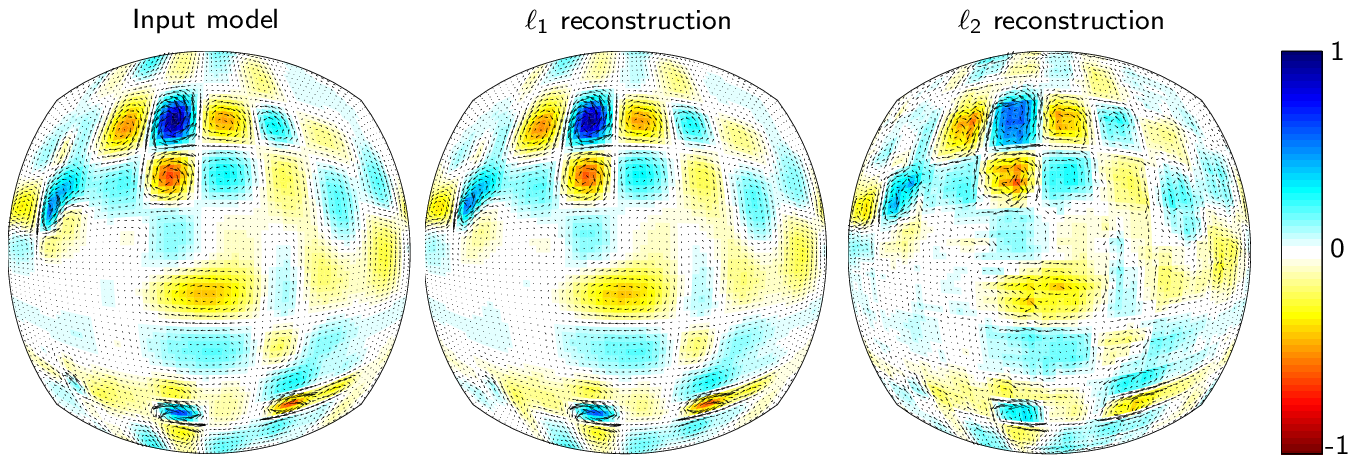}}
\caption{The current density model $\vec{J}$ used in the synthetic sparse
magnetic tomography experiment. Arrows represent the field $\vec{J}$
whereas color indicates the field $F$. Left: the input model. Center:
the $\ell_1$ penalized reconstruction. Right: the $\ell_2$
penalized reconstruction. The $\ell_1$ reconstruction is less noisy and its amplitude is
more faithful to the input model than the $\ell_2$ reconstruction.}\label{megfigure}
\end{figure}

In our simulation we attempt to reconstruct the input model
using an $\ell_1$ penalization approach (to exploit the
sparsity of the model in the wavelet basis), and a traditional
$\ell_2$ penalty approach for comparison. For the former we use
the so-called `Fast iterative soft-thresholding algorithm'
\cite{BeT:2009}, for the latter we use the conjugate gradient
method (the variational equations are linear in that case). The
$\ell_1$ reconstruction of the field $F(\xi,\eta)$ has
a relative reconstruction error of $31\%$, the generic
$\ell_2$ method of $87\%$. Comparing this to the results
on figure \ref{fig:noisetest} we see that, for
$\delta=1000/4096=0.2441$ and $\rho=60/1000=0.06$, the
predicted relative reconstruction error also lies around 30\%
for the $\ell_1$ penalty method. A similar experiment for the
$\ell_2$ penalty method predicts an error around 95\%. The
input model and the two reconstructions are shown in
figure~\ref{megfigure}.

\section{Conclusion}

Sparse recovery can be a very powerful tool for inverse
problems especially when (structured) random matrices or other
matrices that satisfy the Restricted Isometry Property (RIP)
are concerned. However these cases constitute a theoretical
best-case scenario. In this note we investigate the reconstruction
performance when such conditions are loosened and this
`compressed sensing' framework is abandoned. Indeed, most
matrices discussed here do not satisfy the RIP for submatrices
of non-negligible size, but reasonably good reconstruction
results can still be obtained.

By means of extensive numerical simulations we assessed the
performance of the $\ell_1$ recovery method (\ref{l1pen}) for
realistic, non-idealized, linear systems with noisy data.
Various levels of sparsity of the model, different numbers of
data, and different noise levels and spectral behaviors of the
measurement matrices were studied. It was shown that the mean
relative reconstruction error $e$ grows sharply when the
framework of random matrices is abandoned. This means that
accurate reconstruction results can only be achieved for
problems with a more sparse solution. Nonetheless the known
(RIP) bounds from \cite{Rau:2010} seem too restrictive in
practice. These bounds guarantee an accurate reconstruction for
\emph{all} vectors that are sufficiently sparse. But,
although loosening this condition will sometimes provide bad
reconstructions, we showed that it is often possible to
obtain reconstructions of acceptable quality of less sparse
vectors.

Figures~\ref{fig:recon} and \ref{fig:noisetest} allow the
reader to estimate the relative reconstruction error based on
the shape of the measurement matrix (number of rows and
columns), on the expected sparsity of the model and on the
spectral properties of the measurement matrix. These predicted
errors unfortunately are not an upper bound but a lower bound,
in the sense that the simulations did not take into account all
factors that could potentially negatively influence $\ell_1$
penalized reconstructions. Other factors that play an important
role on the quality of reconstructions are e.g. the
distribution of the matrix column norms (i.e. whether basis
vectors lie in the null space of the matrix, in contradiction
to the RIP), and the presence of (nearly) identical columns. In
the latter case, the $\ell_1$ penalty will not impose a unique
minimizer to functional (\ref{l1pen}) and such a penalty is
therefore from the outset not a suitable regularization method
for that inverse problem.

We also conclude that the $\ell_1$ method responds well to an
increased amount of noise: When the noise level increases,
less sparse vectors can be recovered to the same
relative accuracy $e_{2\epsilon}$. The magnetic tomography
example shows the practical usefulness of the numerical study,
and there is good agreement between that toy problem and the
numerical simulations.

\ack This research was supported by VUB GOA-062 and by the
FWO-Vlaanderen grant G.0564.09N. The authors thank Massimo
Fornasier and Francesca Pitolli for discussion
concerning the magnetic tomography problem.

\section*{References}

\bibliography{noisycs}{}
\bibliographystyle{unsrt}

\end{document}